\documentclass[11pt]{article}
\usepackage{amsfonts}
\usepackage{amssymb,amsmath,colordvi,cite}
\usepackage{graphicx,pstricks,psfrag}
\usepackage{CJK}
\newtheorem{theo}{Theorem}[section]

\newtheorem{exam}[theo]{Example}
\newtheorem{lem} [theo]{Lemma}
\newtheorem{coro}[theo]{Corollary}
\newtheorem{prop}[theo]{Proposition}

\makeatletter \@addtoreset{equation}{section}
\@addtoreset{theo}{}\makeatother

\makeindex 

\def\pf{\noindent {\it Proof.} }
\def\qed{\hfill \rule{4pt}{7pt}}

\def\L{\mathop{\rm L}}
\def\D{\mathop{\rm D}}
\begin{document}

\title{\textbf{Abel's Lemma and Identities on Harmonic Numbers}}
\author{Hai-Tao Jin$^1$ and Daniel K. Du$^2$\\[15pt]
$^1$School of Science,\\[5pt]
Tianjin University of Technology and Education, \\[5pt] 
Tianjin 300222, P. R. China \\[5pt]
$^2$Center for Applied Mathematics \\[5pt]
Tianjin University, Tianjin 300072, P. R. China\\[5pt]
$^1$jinht1006@tute.edu.cn, $^2$daniel@tju.edu.cn }
\date{}
\maketitle

\begin{abstract}
Recently, Chen, Hou and Jin used both Abel's lemma on summation by parts and Zeilberger's algorithm
to generate recurrence relations for definite summations. Meanwhile, they proposed the Abel-Gosper method to evaluate some indefinite sums
involving harmonic numbers. In this paper, we use the Abel-Gosper method to prove an identity involving the generalized harmonic numbers. Special cases of
this result reduce to many famous identities. In addition, we use both Abel's lemma and the WZ method to verify and to discover identities involving harmonic numbers. Many interesting examples are also presented.
\end{abstract}

{\noindent\it Keywords}\/: harmonic number, Abel's lemma,  Able-Gosper method, Abel-WZ method

{\noindent\it AMS Classification}\/: 05A19, 33F10, 11B99

\section{Introduction}
The objective of this paper is to employ Abel's lemma on summation by parts and hypergeometric summation algorithms to verify and to discover identities
on the harmonic as well as generalized harmonic numbers.

Recall that for a positive integer $n$ and an integer $r$,
the generalized harmonic numbers in power $r$  are given by
\[
H_n^{(r)}=\sum_{k=1}^{n}\frac{1}{k^r}.
\]
For convenience, we set $H_n^{(r)}=0$ for $n\leq 0$. As usual, $H_n =H_n^{(1)}$ are the classical harmonic numbers. We also define (see \cite{Chu2005})
\[
H_n(x)=\sum_{k=1}^{n}\frac{1}{k+x}, \quad x\neq -1,-2,\ldots
\]
for $n\geq 1$ and $H_n(x)=0$ when $n\leq 0$.
Identities involving these numbers have been extensively
studied and applied in the literature, see, for example, \cite{chosri2011, Chu2005,GKP1994,Paule-Schneider2003,Spies1990}.
Also recall that Abel's lemma \cite{Abel} on summation by parts is stated as follows.

\begin{lem}[Abel's lemma]
For two arbitrary sequences
$\{a_k\}$ and $\{b_k\}$, we have
\[
\sum_{k=m}^{n-1} (a_{k+1}-a_k)b_k =\sum_{k=m}^{n-1}a_{k+1}(b_k-b_{k+1})+a_{n}b_{n}-a_{m}b_{m}.
\]
\end{lem}
For a sequence $\{\tau_k\}$, define the forward difference operator $\Delta$ by
\[
\Delta \tau_k = \tau_{k+1} - \tau_k.
\]
Then Abel's lemma can be written as
\begin{equation}\label{Abel'lemma}
\sum_{k=m}^{n-1}b_k\Delta a_k=-\sum_{k=m}^{n-1}a_{k+1}\Delta b_k+a_{n}b_{n}-a_{m}b_{m}.
\end{equation}

Graham, Knuth and Patashnik \cite{GKP1994} reformulated Abel's lemma in terms of finite
calculus to evaluate several sums on harmonic numbers. Recently, Chen, Hou and Jin \cite{Chen-Hou-Jin}
proposed \emph{the Abel-Gosper method} and derived some identities on harmonic numbers.
The idea can be explained as follows. Let $f_k$ be a hypergeometric term, i.e., $f_{k+1}/f_k$ is a rational function of $k$.
First, we use Gosper's algorithm \cite{PWZ1996} to find a hypergeometric term $a_k$ (if it exists) satisfying $\Delta a_k = f_k$.
Then, by Abel's lemma, we have
\begin{equation}\label{H-trans}
\sum_{k=m}^{n-1} f_k H_k = \sum_{k=m}^{n-1} H_k \Delta a_k = -\sum_{k=m}^{n-1} \frac{a_{k+1}}{k+1}
+ a_n H_n - a_m H_m.
\end{equation}
Hence we can transform a summation involving harmonic numbers into a hypergeometric summation. For example, let
\[
S(n) = \sum_{k=1}^{n} H_k.
\]
We have
\begin{equation*}
S(n) = \sum_{k=1}^n H_k \Delta k = - \sum_{k=1}^n (k+1) \Delta H_k + (n+1)H_{n+1} - H_1 = (n+1) H_n - n.
\end{equation*}
In this framework, they combine both Abel's lemma and Zeilberger's algorithm to find recurrence relations
for definite summations involving non-hypergeometric terms. For example, they can prove the Paule-Schneider identity \cite{Paule-Schneider2003}
\[
\sum_{k=0}^n\left(1+3(n-2k)H_k\right){n\choose k}^3 = (-1)^n,
\]
and Calkin's identity \cite{Calkin1994}
\[
\sum_{k=0}^n\left(\sum_{j=0}^k{n\choose j}\right)^3 = n 2^{3n-1} +2^{3n} -3n2^{n-2}{2n\choose n}.
\]

In this paper, we use the Abel-Gosper method to generalize the following
well-known inversion formula (see, for example \cite[(1.46)]{Gould1972})
\begin{equation}\label{Binom1-eq-1}
\sum_{k} (-1)^{k-1}\binom{n}{k} H_k=\frac{1}{n}.
\end{equation}
To be specific, we have
\begin{theo}\label{Main-1}
Let $m,s,p,n$ are nonnegative integers with $n\geq p$ and $m\geq 1$, then
\begin{equation}\label{Main-1-eq}
\sum_{k=p}^{n} (-1)^{k-1}\binom{n}{k} \binom{k}{p} H_{mk+s}(x)= 
 \left\{
\begin{array}{cc}
\frac{(-1)^{p}m^{n-p-1}n!}{(n-p)p!}\sum\limits_{i=1}^{m}\frac{1}{\prod_{u=p}^{n-1}(mu+s+x+i)},&
n>p,\\[7pt]
(-1)^{p-1}H_{mp+s}(x),& n=p.
\end{array}
\right.
\end{equation}
\end{theo}
It is readily to see that identity \eqref{Main-1-eq} reduce to inversion formula \eqref{Binom1-eq-1} by setting  $p=0, m=1, s=0$ and $x=0$.
More interesting special cases of \eqref{Main-1-eq} can be found in Section~2.

In addition, by combining Abel's lemma with the WZ method, we establish the \emph{Abel-WZ method} to construct
identities on harmonic numbers from known hypergeometric identities. For example, we shall reestablish the following identity due to Prodinger \cite{Prodinger2008}.
\begin{equation*}\label{Prodinger2008}
\sum_{k=0}^{n} (-1)^{n-k}\binom{n}{k}\binom{n+k}{k}  H_{k}^{(2)}=2\sum_{k=1}^{n}\frac{(-1)^{k-1}}{k^2}.
\end{equation*}

The paper is organized as follows.
In Section~2, we shall give a proof of Theorem \ref{Main-1} by the Abel-Gosper method. Special cases of Theorem \ref{Main-1} and more examples are
also displayed. In Section~3, we introduce the Abel-WZ method and then construct many interesting identities on harmonic numbers from hypergeometric identities.


\section{The Abel-Gosper Method}\label{Sec-Abel-Gosper}
We first make use of the Abel-Gosper method to prove Theorem \ref{Main-1}.

{\noindent {\it Proof of Theorem \ref{Main-1}.}}
Let
\[
S_{m,s,p}(n,x)=\sum_{k=p}^{n} (-1)^{k-1}\binom{n}{k} \binom{k}{p} H_{mk+s}(x),
\]
and
\[
F(n,k)=(-1)^{k-1}\binom{n}{k} \binom{k}{p}.
\]
By Gosper's algorithm, we have
\[
F(n,k)=\Delta_k G(n,k),
\]
where
\[
G(n,k)=\frac{(-1)^k(k-p)}{n-p}\binom{n}{k}\binom{k}{p}.
\]
Thus it follows that
\[
S_{m,s,p}(n,x)=\sum_{k} \Delta_k G(n,k) H_{mk+s}(x).
\]
Employing Abel's lemma and noticing the boundary values, we find that
\[
S_{m,s,p}(n,x)=-\frac{1}{n-p}\sum_k
(-1)^{k-1}(k+1-p)\binom{n}{k+1}\binom{k+1}{p}\sum_{i=1}^{m}\frac{1}{mk+s+i+x}.
\]
For $1\leq i \leq m$, set
\[
S_i(n)=\sum_k
(-1)^{k-1}(k+1-p)\binom{n}{k+1}\binom{k+1}{p}\frac{1}{mk+s+i+x}.
\]
Then Zeilberger's algorithm (see \cite{PWZ1996}) returns the
recurrence equation
\[
(mn+s+i+x)S_i(n+1)-m(n+1)S_i(n)=0.
\]
By the initial value
\[
S_i(p+1)=\frac{(-1)^{p+1}(p+1)}{mp+s+i+x},
\]
we obtain
\[
S_i(n)=(-1)^{p+1}\frac{m^{n-p-1}n!}{p!\prod_{u=p}^{n-1}(mu+s+i+x)}.
\]
Equation \eqref{Main-1} is then established by noticing that
\[
S_{m,s,p}(n,x)=-\frac{1}{n-p}\sum_{i=1}^{m}S_i(n),\quad n>p
\]
and $S_{m,s,p}(p)=(-1)^{p-1}H_{mp+s}(x)$. \qed

Now let us show some special cases of Theorem \ref{Main-1}.
By setting $m=1$ and $x=0$, \eqref{Main-1-eq} reduces to the following identity.
\begin{coro}\label{Main-1-sp1}
Let $n,p$, $s$ be nonnegative integers and $n>p$, then we have
\begin{equation}\label{Main-1-sp1-eq}
\sum_{k=p}^{n}(-1)^{k-1}\binom{n}{k}\binom{k}{p}H_{k+s}=\frac{(-1)^p \binom{p+s}{s}}{(n-p)\binom{n+s}{s}}.
\end{equation}
\end{coro}
The special cases $p=0$ and $s=0$ of \eqref{Main-1-sp1-eq} are given in \cite{Spies1990,WangWeiPing2010}.

By setting $m=2,s=0$ and $x=0$ in \eqref{Main-1-eq}, we are led to the following identity .
\begin{coro}\label{Main-1-sp2}
Let $n, p$ be nonnegative integers and $n>p$, then we have
\begin{equation}\label{Main-1-sp2-eq}
\sum_{k=p}^{n} (-1)^{k-1}\binom{n}{k} \binom{k}{p}
H_{2k}=\frac{(-1)^p}{(n-p)}\left(\frac{1}{2}+\frac{
2^{2n-2p-2}\binom{2p}{p}}{\binom{2n-1}{n-1}}\right).
\end{equation}
\end{coro}
Using the relation $k^2=2\binom{k}{2}+\binom{k}{1}$ and the cases $p=1,2$ of \eqref{Main-1-sp2-eq}, we arrive at an identity due to Sofo \cite{Sofo2009}.
\begin{equation}\label{Main-1-sp2-eq1}
\sum_{k} (-1)^{k-1}\binom{n}{k} k^2
H_{2k}=\frac{n}{2(n-1)(n-2)}+\frac{2^{2n-4}}{(n+2)\binom{2n-1}{n-3}}, \quad n>2.
\end{equation}

Note that we can also derive identities involving the generalized harmonic numbers $H_n^{(r)}$ from Theorem \ref{Main-1}. To this end, we need the operators
$\L$ and $\D$ which are defined by $\L f(x)=f(0)$ and $\D f(x)=f'(x)$. It is readily to see that
\[
\L {\D}^m H_n(x)=(-1)^{m}m!H^{(m+1)}_n.
\]

By setting $m=1$ and $p=0$ in \eqref{Main-1-eq}, we get the following result (see \cite{Larcombe2005}).
\begin{coro}
\begin{equation}\label{Main-2-sp-eq}
\sum_{k=0}^{n} (-1)^{k-1}\binom{n}{k}H_{k+s}(x)=\frac{n!}{n
(s+x+1)_n}.
\end{equation}
\end{coro}
Then applying the operator $\L \D$ to both sides of \eqref{Main-2-sp-eq},
we obtain a formula given in \cite{WangWeiPing2010}
\begin{equation}\label{Main-2-sp-eq-1}
\sum_{k=0}^{n}(-1)^{k-1}\binom{n}{k}H^{(2)}_{k+s}=
-\frac{1}{n}(H_s-H_{n+s})\binom{n+s}{s}^{-1}.
\end{equation}
Furthermore, applying the operator $\L \D^2$ to both sides of \eqref{Main-2-sp-eq}
gives
\[
\sum_{k=0}^{n}(-1)^{k-1}\binom{n}{k}H^{(3)}_{k+s}=
\frac{1}{2n}\Big((H_{n+s}-H_s)^2+H^{(2)}_{n+s}-H^{(2)}_s\Big)\binom{n+s}{s}^{-1}.
\]

More generally, \eqref{Main-2-sp-eq} leads to the following inversion formula by applying the operator $\L \D^m$ to its both sides.
\begin{prop}
For positive integers $n$ and $m$, we have
\begin{equation}\label{Main-2-sp-eq-2}
\sum_{k=1}^{n}(-1)^{k-1}\binom{n}{k}H^{(m+1)}_{k}=\frac{1}{n}\sum_{1\leq
j_1\leq j_2\leq\cdots\leq j_m\leq n}\frac{1}{j_1 j_2\cdots j_m}.
\end{equation}
\end{prop}
\pf Setting $s=0$ in \eqref{Main-2-sp-eq} and applying the operator $\L {\D}^m$
to its both sides, we have
\[
(-1)^{m}m! \sum_{k} (-1)^{k-1}\binom{n}{k}
H^{(m+1)}_k=\frac{n!}{n}LD^m\frac{1}{(x+1)_n}.
\]
By the partial fraction decomposition
\[
\frac{1}{(x+1)_n}=\sum_{k=1}^{n}\frac{1}{(x+k)\prod\limits_{1\leq
j\neq k\leq n}(j-k)},
\]
we find
\[
\frac{n!}{n}\L {\D}^m\frac{1}{(x+1)_n}=(-1)^m
m!\frac{1}{n}\sum_{k=1}^{n}\binom{n}{k}\frac{(-1)^{k-1}}{k^m}.
\]
Finally, using Dilcher's formula \cite{Dilcher1995}
\[
\sum_{k=1}^{n}\binom{n}{k}\frac{(-1)^{k-1}}{k^m}=\sum_{1\leq j_1\leq
j_2\leq\cdots\leq j_m\leq n}\frac{1}{j_1 j_2\cdots j_m},
\]
we arrive at \eqref{Main-2-sp-eq-2}. \qed

Similarly, we can use the Abel-Gosper method to find many other identities.
Here are some examples.

\begin{exam}
For $n\in \mathbb{N}$ and $x\in \mathbb{C}\setminus \{-1,-2,\ldots\}$, we have
\begin{align*}
\sum_{k=0}^{n}\frac{(x+1)_k}{k!}H_k&=\frac{1}{x+1}\left(1+\frac{(x+1)_{n+1}}{n!}\left(H_n-\frac{1}{x+1}\right)\right),\\[6pt]
\sum_{k=0}^{n}\frac{k!}{(x+1)_{k}}H_k&=\left\{
\begin{array}{ll}
\frac{1}{(x-1)^2}\left(x-\frac{n!}{(x+1)_{n}}\left((x-1)(n+1)H_{n}+n+x\right)\right), &\mbox{if}\quad x\neq 1,\\[6pt]
\frac{H^2_{n+1}-H^{(2)}_{n+1}}{2},&\mbox{if}\quad x= 1.
\end{array}
\right.
\end{align*}
\end{exam}
We remark that the second identity  also holds when $x$ is a negative integer. In this case, it it equivalent to
the following formula (see \cite[Exercise 6.53]{GKP1994}).
\[
\sum_{k=0}^{n}\frac{(-1)^k}{\binom{m}{k}}H_{k}=
\frac{(-1)^n}{\binom{m}{n}}\left[\frac{n+1}{m+2}H_{n}+\frac{m+1-n}{(m+2)^2}\right]-\frac{m+1}{(m+2)^2},
\]
where $m,n\in \mathbb{N}$ and $n\leq m$.

\begin{exam}
For $n,m,p\in \mathbb{N}$, we have the following three identities.
\begin{align*}
&\sum_{k=0}^{n}(-1)^{k-1}\frac{\binom{n}{k}}{\binom{k+p}{p}}H_{k+p}=\frac{n-p(n+p)H_p}{(n+p)^2},\\[6pt]
&\sum_{k=0}^{n}(-1)^{k-1}\frac{k \binom{n}{k}}{\binom{k+p}{p}}H_{k}=\frac{pn(1+H_{p-1}-H_{n+p-2})}{(n+p)(n+p-1)},\quad p\geq 2,\\[6pt]
&\sum_{k=0}^{n}(-1)^{k-1}\frac{k^2 \binom{n}{k}}{\binom{k+p}{p}}H_{k}=\frac{pn((n-p)(H_{n+p-3}-H_{p-1})-(2n-p))}{(n+p)(n+p-1)(n+p-2)},\quad p\geq 3.
\end{align*}
\end{exam}
We remark that the first formula is due to Sofo \cite{Sofo2008} and the remaining two are obtained by Chu \cite{Chu2012}.

Using the Abel-Gosper method iteratively, we can prove the following identity.
\begin{exam}
For $n,p\in \mathbb{N}$ and $n>p$, we have
\begin{equation*}\label{new-2}
\sum_{k} (-1)^{k-1}\binom{n}{k} \binom{k}{p}H_k^{2}=\frac{(-1)^p}{n-p}(H_n-2H_{n-p-1}+H_p).
\end{equation*}
\end{exam}

\section{The Abel-WZ Method}

In this section, we shall illustrate how to combine Abel's lemma with the WZ method to derive
identities on harmonic numbers.

Recall that a pair of hypergeometric
functions $(F(n,k),G(n,k))$ is called a WZ pair if the following  WZ equation holds
\[
F(n+1,k)-F(n,k)=G(n,k+1)-G(n,k).
\]
For a given $F(n,k)$, the WZ method will give such $G(n,k)$ if it exists, see for example \cite{PWZ1996}.
Now we are ready to describe the Abel-WZ method. Assuming that we have the following hypergeometric identity.
\[
\sum_k F(n,k)=f(n).
\]
In most cases, we can obtain a WZ pair
\[
\left(\frac{F(n,k)}{f(n)},G(n,k)\right).
\]
Then for the sum $S(n)=\sum_{k\geq 0} F(n,k)b_k$, where $b_k$ is harmonic number, we have
\[
\frac{S(n+1)}{f(n+1)}-\frac{S(n)}{f(n)}=\sum_k (G(n,k+1)-G(n,k))b_k.
\]
Denote by $U(n)=\sum_k (G(n,k+1)-G(n,k))b_k$. Then by Abel's lemma, we have (here we omit the boundary values)
\[
U(n)=-\sum_k G(n,k+1)\Delta_k b_k.
\]
Again, if $\Delta_k b_k$ is hypergeometric, $U(n)$ can be treated by Zeilberger's algorithm.
Moreover, if $U(n)$ can be expressed in closed form, we then establish an identity of the form
\[
S(n)=f(n)\sum_{k\leq n-1} U(k).
\]

We begin by an identity due to Prodinger \cite{Prodinger2008}.
\begin{exam}
For $n\in \mathbb{N}$, we have
\begin{equation}\label{Prodinger2008}
\sum_{k=0}^{n} (-1)^{n-k}\binom{n}{k}\binom{n+k}{k}  H_{k}^{(2)}=2\sum_{k=1}^{n}\frac{(-1)^{k-1}}{k^2}.
\end{equation}
\end{exam}
\pf Denote the left side of \eqref{Prodinger2008} by $S(n)$. For $F(n,k)=(-1)^{n-k}\binom{n}{k}\binom{n+k}{k}$, the WZ method gives
\[
F(n+1,k)-F(n,k)=G(n,k+1)-G(n,k),
\]
where
\[
G(n,k)=\frac{2(-1)^{n-k}k^2{n\choose k}{n+k\choose k}}{(n-k+1)(n+1)}.
\]
Multiplying both sides of the WZ equation by $H_{k}^{(2)}$ and summing over $k$ gives
\[
S(n+1)-S(n)=\sum_{k}(G(n,k+1)-G(n,k))H_{k}^{(2)}.
\]
Then applying Abel's lemma to the right hand side of the above identity and noting the boundary values, we have
\begin{align*}
S(n+1)-S(n)&=\sum_k \frac{-G(n,k+1)}{(k+1)^2}\\
           &=\sum_{k\geq 0} \big( T(k+1)-T(k) \big) \\[5pt]
           &=-T(0)=2\frac{(-1)^n}{(n+1)^2},
\end{align*}
where
\[
T(k)=\frac{2(-1)^{n-k-1}(k+1)^2{n\choose k+1}{n+k+1\choose k+1}}{(n-k)(n+1)^3}.
\]
Thus we have
\[
S(n)=S(0)+2\sum_{k=1}^{n}\frac{(-1)^{k-1}}{k^2}.
\]
By the initial value $S(0)=0$, we complete the proof.\qed

The underlying hypergeometric identity of the above theorem is the special case $p=0$ of
\[
\sum_{k=0}^{n}(-1)^k\binom{n}{k}\binom{n+k}{k}\binom{k}{p}=(-1)^n\binom{n+p}{p}\binom{n}{p},
\]
which enables us to establish the following identities.

\begin{exam}
For $n,p\in \mathbb{N}$ and $n\geq p$, we have
\begin{align*}
\sum_{k=0}^{n}(-1)^{n-k}{n\choose k}{n+k\choose k} H_{2k}&=3H_{n}-H_{\lfloor \frac{n}{2}\rfloor}, \\[6pt]
\sum_{k=0}^{n}(-1)^k\binom{n}{k}\binom{n+k}{k}\binom{k}{p}H_k&=(-1)^n\binom{n+p}{p}\binom{n}{p}(2H_n-H_p),\\[6pt]
\sum_{k=0}^{n}(-1)^k\binom{n}{k}\binom{n+k}{k}\binom{k}{p}H_{n+k}&=(-1)^n\binom{n+p}{p}\binom{n}{p}(H_{n+p}+H_n-H_p).
\end{align*}
The cases $p=0,1$ of the last two formulas can be found in \cite{Prodinger2008} and \cite{Osburn} respectively.
\end{exam}

The following contents of this section consist of several selected examples.

\begin{exam}
From the binomial theorem $\sum_{k}\binom{n}{k}\lambda^{n-k}\mu^k=(\lambda+\mu)^n$, we can derive
the following formula due to Boyadzhiev \cite{Boyadzhiev2009}.
\begin{equation*}\label{Boyadzhiev}
\sum_{k=1}^{n}\binom{n}{k}H_k \lambda^{n-k}\mu^k=(\lambda+\mu)^n H_n-
\left(\lambda (\lambda+\mu)^{n-1}+\frac{\lambda^2}{2}(\lambda+\mu)^{n-2}+\cdots+\frac{\lambda^n}{n} \right).
\end{equation*}
\end{exam}

\begin{exam}
From identity
\[
\sum_{k=p}^{n}{n\choose k}^2{k\choose p}={2n-p\choose n}{n\choose p},
\]
we can derive
\begin{equation*}\label{eq-wz-hk}
\sum_{k=p}^{n}{n\choose k}^2{k\choose p}H_k={2n-p\choose n}{n\choose p}(2H_n-H_{2n-p}).
\end{equation*}
The special cases $p=0$ and $p=1$ are due to Paule and Schneider \cite{Paule-Schneider2003}.
\end{exam}

\begin{exam}
From the identities
\[
\sum_{k=0}^{2n}(-1)^k{2n\choose k}^2=(-1)^n{2n\choose n}
\]
and
\[
\sum_{k=0}^{2n}(-1)^k{2n\choose k}^3=(-1)^n\frac{(3n)!}{n!^3},
\]
we have
\begin{eqnarray*}
\sum_{k=0}^{2n}(-1)^k{2n\choose k}^2 H_k&=&(-1)^n{2n\choose n}\frac{H_n+H_{2n}}{2},\label{eq-wz-41}\\ [6pt]
\sum_{k=0}^{2n}(-1)^k{2n\choose k}^3 H_k&=&(-1)^n\frac{(3n)!}{n!^3}\frac{H_n+2H_{2n}-H_{3n}}{2},\label{eq-wz-42}\\ [6pt]
\sum_{k=0}^{2n}(-1)^k{2n\choose k}^3 H_{k}^{(2)}&=&(-1)^n\frac{(3n)!}{n!^3}\frac{H_{n}^{(2)}+H_{2n}^{(2)}}{2}.\label{eq-wz-43}
\end{eqnarray*}
The last two formulas can be found in \cite{Driver-Schneider20062} and \cite{ChuFu2009} respectively.
\end{exam}

\vspace{.2cm} \noindent{\bf Acknowledgments.}
We wish to thank Qing-Hu Hou for helpful comments and discussions. This work was
supported by the National Science Foundation of China (Tianyuan Fund for Mathematics) and
the Project Sponsored by the Scientific Research Foundation of Tianjin University of Technology and Education.


\end{document}